\documentclass{article}

\def\eps{\varepsilon}
\def \Z {{\mathbf {Z}}}

\def \N {{\mathbf {N}}}

\usepackage[T2A]{fontenc}
\usepackage[cp1251]{inputenc}
\usepackage[tbtags]{amsmath}
\usepackage{amsfonts,amssymb}

\begin{document}

\begin{center}   
{ \LARGE \bf
Lecture notes on ergodic transformations
\vspace{5mm}

Valery V. Ryzhikov
}
\end{center} 
 \rm

\vspace{5mm}
Our aim is to introduce the reader to some classical properties of ergodic transformations and some problems (many of them are not so simple). The text is based on lectures that the author gave at the Faculty of Mechanics and Mathematics of Moscow State University for senior students wishing to become familiar with the theory of ergodic transformations. Several introductory lectures were held at the international summer school at Moscow State University on July 23-25, 2024, what was the reason for publishing this version.  

\Large
\vspace{15mm}
\ \ \ \ \ \ \ \ \ \ \ \ \ \ \ \ \ \ \ \ \ \
  \ \ \ \ \ \ \ \ \ \ \ \ \ \ \ \ \ \ {\bf \LARGE Content} 

 \ \ \ \ \ \ \ \ \ \ \ \ \ \ \ \ \S 1. \ Introduction.
 
 \ \ \ \ \ \ \ \ \ \ \ \ \ \ \ \ \S 2. \ {Properties of transformations equivalent to ergodicity.}

 \ \ \ \ \ \ \ \ \ \ \ \ \ \ \ \ \S 3.  \ {Birkhoff's Theorem.}

 \ \ \ \ \ \ \ \ \ \ \ \ \ \ \ \  \S 4. \ Properties equivalent to weak mixing.

 \ \ \ \ \ \ \ \ \ \ \ \ \ \ \ \  \S 5. \ On typical properties of transformations.

 \ \ \ \ \ \ \  \ \ \ \ \ \ \ \ \  \S 6. \ Lego to constuct transformations. 

 \ \ \ \ \ \ \ \ \ \ \ \ \ \ \ \ \S 7. \ Typical entropy invariants.

 \ \ \ \ \ \ \ \ \ \ \ \ \ \ \ \ \S 8. \ {Poisson suspensions with completely positive $P$-entropy.}  

 \ \ \ \ \ \ \ \ \ \ \ \ \ \ \ \ \S 9. \ Spectral theorem for unitary operators.

 \ \ \ \ \ \ \ \ \ \ \ \ \ \ \ \ \S 10. Compact factors, Kronecker algebra.

 \ \ \ \ \ \ \ \ \ \ \ \ \ \ \ \ \S 11. Progression recurrence for weakly mixing transformations.

 \ \ \ \ \ \ \ \ \ \ \ \ \ \ \ \  \S 12.  Double recurrence for ergodic  transformations.

\Large
\section{Introduction}
For brevity, by \it transformation \rm  we mean a measure-preserving invertible  transformation $T$ of the Lebesgue space $(X,{\mathcal B }, \mu)$. Such transformations are often called automorphisms of a measure space. Now we will discuss the standard probability space, i.e. $\mu(X)=1$. Transformations that differ on a set of measure 0 are identified.

\bf Ergodicity. \rm Transformation $T$ is ergodic if every measurable set $A$ such that $A=TA$,
has measure 0 or 1. In other words, the phase space $X$ is indecomposable: cannot be "cut"  into two nontrivial invariant parts.

\bf Weak mixing. \rm Transformation $T$ is called weakly mixing if
there is a sequence $n_i\to\infty$ such that
$$\forall A,B\in{\mathcal B} \ \ \ \mu(T^{n_i}A\cap B)\to \mu(A)\mu(B), \ i\to\infty.$$

\bf Mixing. \rm Transformation $T$ is called mixing if
$$\forall A,B\in{\mathcal B} \ \ \ \mu(T^iA\cap B)\to \mu(A)\mu( B), \ i\to\infty.$$
(This property is not typical in the Baire sense, we will prove this later.)

\vspace{2mm}
\bf Example of ergodic transformation. \rm Let $a$ be a real number, consider the shift
$T_a\ : [0,1)\to [0,1)$
$$T_a(x)=\{x+a\}:= x+a \ (mod\ 1).$$ \it The transformation $T_a$ is ergodic with respect to Lebesgue measure on $[0,1)$
iff $a$ is irrational. \rm Hint: find an elementary proof of this fact based on the fact that any (!) measurable
set $A$ of positive measure intersects some segment of measure ${\varepsilon}$ to a measure greater than $0.9 {\varepsilon}$. If $TA=A$, then  can get   $\mu(A)>0.8$, hence, $\mu(A)=1$ or $\mu(X\setminus A)>0.8$
(the latter is impossible since $\mu(A)>0.8$).

Let's consider a more general example: a group shift on a torus
$$ T_{(a_1,\dots, a_n)}(x_1,\dots,x_n)=(\{x_1+a_1\},\dots, \{x_n+a_n\}),$$
where $\{x+a\}$ denotes the fractional part of $x+a$.
What conditions should be imposed on $a_1,\dots, a_n$ so that the shift is ergodic with respect to the Lebesgue measure on the torus?

\vspace{2mm}
\bf An example of a mixing transformation. \rm Consider the Bernoulli shift $T: {\mathbf Z}_2^{{\mathbf Z}}\to {\mathbf Z}_2^{{\mathbf Z}}$ on the space of
two-sided sequences taking the values from $\{0, 1\}$. For $x$  such a sequence, we consider the shift
$$ T(x)_z=x_{z-1}.$$
The measure $\mu$ on $X={{\mathbf Z}}_2^{{\mathbf Z}}$ is defined as follows: the set of all sequences with
fixed values for on a  finite set $F$ of indices $z$ has measure $\frac 1{2^{|F|}}$.
It is natural to call such sets cylinders; they form a semiring on which we have defined a countably additive measure called Bernoulli mesure of type $\left(\frac 1 2, \frac 1 2\right)$.
It is sufficient to verify the mixing property on cylinders. If $A$ and $B$ are cylinders, then
$$\mu(T^iA\cap B)=\mu(A)\mu( B)$$
for all sufficiently large $i$, which follows easily from the definition of the measure.

Note that the shift $T$ is an automorphism of the compact group ${{\mathbf Z}}_2^{{\mathbf Z}}$ that preserves its Haar measure (here, this is a probability measure invariant with respect to all group shifts). And this  Bernoulli-Haar measure is also  Lebesgue on a square.  Indeed, the automorhism $T$ is also called the baker's transformation of $X$ (as a square), and $T$ maps a point  to another point of the square as follows:
$$ T\left(\sum_{i=1}^{\infty} \frac {x_i}{2^i},\ \sum_{i=1}^{\infty} \frac {x_{-i+1}}{2^i}\right)=\
\left(\sum_{i=1}^{\infty} \frac {x_{i-1}}{2^i},\ \sum_{i=1}^{\infty} \frac {x_{-i}}{2^i}\right).$$
The baker's transformation $T:[0,1]^2\to [0,1]^2$ preserves Lebesgue measure on a square $[0,1]^2$. Locally, it contracts  along one coordinate and expands  along the other.

The property of weak mixing follows directly from the property of mixing, which implies  ergodicity.
Indeed, if $TA=A$, then $T^nA=A$, then taking into account weak mixing for $B=A$, we have
$$\mu(A)=\mu(T^{n_i}A\cap A)\to \mu(A)\mu(A), \ \ \ \mu(A)\in \{0,1\}.$$
This means that $T$ is ergodic (the invariant set is all or nothing).

Note that circle rotations (shifts on tori) do not have weak mixing.
Weak mixing implies the following property: the sets $T^{n_i}A$ asymptotically have a nonempty intersection with any
fixed sets of positive measure. However, a small segment under rotation cannot intersect two distant segments simultaneously --
the absence of the mixing  in this case is obvious.

\vspace{2mm}
\bf Koopman operators. \rm Let us consider for an invertible measure-preserving transformation $T$
a unitary operator $T$ acting in $L_2(X,{\mathcal B},\mu)$
according to the simple formula $$Tf(x)=f(T(x)).$$
We denote the operator and the transformation in the same way, this does not lead to misunderstandings. In what follows,
we write $Tx$ instead of $T(x)$. The mentioned shifts by ${{\mathbf Z}}_2^{{\mathbf Z}}$ and by ${{\mathbf Z}}_3^{{\mathbf Z}}$ are isomorphic as operators:
they are identically structured permutations of suitable orthonormal systems.
If the reader is familiar with the characters of commutative groups (or with Rademacher-Walsh systems), then
he will quickly find these systems and shifts on them. However, as Kolmogorov showed, unitary isomorphism is not induced by any measure-preserving transformation. Two transformations $S$, $T$ are called metrically isomorphic if there is a transformation $R$ such that $S=R^{-1}TR$.

Denote by $\Theta$ the orthoprojection onto the space of constants in $L_2$.
The mixing property is equivalent to the convergence
$$T^i\to_w \Theta\ : \  \ \int_X T^if\, g \,d\mu \ \to\ \int_X f\, d\mu \int_X g\, d\mu, \ i\to\infty,$$
for all $f,g\in L_2(X,\mu)$.
If the operators $T^{n_i}$ converge weakly to $\Theta$:
$$T^{n_i}\to_w \Theta$$ for some sequence ${n_i}$, we call such $T$ weakly mixing.

The ergodicity property can also be formulated in terms of weak limits:
$$\frac 1 N \sum_1^N T^i \to_w \Theta.$$
We will verify  this  later.

\newpage
\bf Infinite family of non-isomorphic weakly mixing non-mixing transformations. \rm
If we are familiar with the ergodic rotation of the circle and the Bernoulli shift, and also know the Rokhlin-Halmos lemma,
then in fact we are able to  construct an uncountable family of non-isomorphic transformations. Let's try.

\vspace{2mm}
\bf Rokhlin-Halmos lemma. \it For an ergodic transformation $T$, ${\varepsilon}>0$
and a natural $H$, there exists a measurable set $B$ such that
the sets $T^iB$ are disjoint for $0\leq i\leq H-1$ and
$$\mu\left({\bigsqcup}_{i=0}^{H-1} T^iB\right)>1-{\varepsilon}.$$\rm


\vspace{2mm}
 \bf Exercises. \rm Find a set $A$ such that for some $n\gg H$ we have
$$\mu\left({\bigsqcup}_{i=0}^{n} T^iA\right)>0.$$

\vspace{2mm} 
  -- Represent the entire space $X$ as a union of disjoint (high) towers:
$$X={\bigsqcup}_{n>N}{\bigsqcup}_{i=0}^{n} T^iA_n, \ \ N\gg H,$$
and  prove the lemma.

\vspace{2mm} 
-- Prove the Lehrer-Weiss  $\eps$-free lemma:

 \it If $T^n$, $n>0$, is ergodic, then for any set $A$, $\mu(A)>0$,  
there exists a measurable set $B$ such that
$$X\setminus A\ \subset \ {\bigsqcup}_{i=0}^{n-1} T^iB.$$\rm

\vspace{2mm}
\bf Corollary of the Rokhlin-Halmos lemma. \it For ergodic transformations $T,T'$ and any $\delta>0$, one can modify the transformation $T$ on a set
of measure less than $\delta$, such that the modified transformation is isomorphic to (conjugate to) $T'$. \rm

\vspace{2mm}
Hint: For $T$ and $T'$ we apply R-H lemma with $\mu(B)=\mu(B')$ (this is easy to ensure)
and ${\varepsilon}$ and $H$ such that ${\varepsilon}+\frac 1 H <<\delta$.
Next we find $R$ such that on the set ${\bigsqcup}_{i=0}^{H-2} {T'}^iB'$ the conjugation $R^{-1}TR$ coincides with $T'$.

\bf Weakly mixing rigid transformations. \rm
An ergodic rotation of the circle $S$ has the rigidity property: $S^{k(j)}\to_w I$, where $I$ is the identity operator, $k(j)\to\infty$
(note that weak convergence of unitary operators to a unitary operator turns out to be strong convergence). The rigidity property, like the mixing property, is invariant under conjugation. We will use this now.

Consider a sequence of transformations $T_p$ such that
$$\sum_p\mu(\{x: T_p(x)\neq T_{p+1}\})<\infty,$$
then $T_p$ converge almost everywhere to the transformation $T$.
We will choose $T_{2p}$ isomorphic to the rigid transformation, and $T_{2p+1}$ isomorphic to the mixing transformation,
Since the rate of convergence of the series is allowed to be chosen arbitrarily fast,
for the limit transformation $T$ we can ensure convergence
$$ T_{2p+1}^{k(2p+1)}\approx_w T^{k(2p+1)}\to_w \Theta, \ \  \ \ \ T_{2p}^{k(2p)}\approx_w T^{k(2p)}\to_w I, \ p\to\infty.$$
The above notation $T_{2p+1}^{k(2p+1)}\approx_w T^{k(2p+1)}$ means that $T_{2p+1}^{k(2p+1)}- T^{k(2p+1)}\to_w 0$.

Thus, by varying the transformation on sets of very small measure, we can control the rigid and mixing sequences and thus obtain a continuum of non-isomorphic weakly mixing but non-mixing transformations.  We recall that for isomorphic transformations, both rigid and mixing sequences coincide.

\vspace{2mm} 

\bf Amusing problems. \rm

\vspace{2mm}
--Prove that the circle rotation and the Bernoulli shift are metrically isomorphic to their inverses, respectively.

\vspace{2mm}
-- The full group $[T]$ of $T$ consists of all so-called generalized powers, that is,  $S\in [T]$ means that the invertible transformation $S$ has the form $S(x)=T^{n(x)}x$,
where $n(x)$ is some measurable function taking integer values.
(The famous Dye's theorem states that the full group of an ergodic transformation contains isomorphic copies of all ergodic transformations!) Prove that  an ergodic  transformation $T$ is not conjugate to its inverse  in the group $[T]$.

\vspace{2mm}
-- Following  Arnaldo Nogueira, we rotate the circle by some angle. Then on some  arc of this circle we make the following  flip:
$\alpha+\varphi\to \alpha-\varphi$, where $\alpha$ is the center of the arc (do nothing outside the arc). As a result of the above composition
we get a transformation of the circle. Prove that this transformation is  periodic (this  surprises   at times).

-- Let's consider now the composition of an ergodic group shift $T$ on the 2-torus and  a  rotation $S$ by 180 degrees of some  small disk in the torus. Outside the disk, the transformation $S$ is identical. Prove that the  composition $ST$ is periodic.

\vspace{2mm}
-- Now we   consider the composition $ST$  of the  baker transformation $T:[0,1]^2\to [0,1]^2$ with the  flip $S:[0,1]^2\to [0,1]^2$, which is the  symmetry of a triangle with vertices $(0,0), (\eps, 0),$ $(0,\eps)$  relative to the diagonal in a square  ($\eps>0$). Outside the triangle $S=Id$. Prove that  the Lyapunov exponents (see Valery Oseledets multiplicative ergodic theorem)  for $ST$ are zero.
If this seems complicated, prove simply that almost all points of the square are periodic with respect to $ST$.

\section{Properties of transformations equivalent to ergodicity}
\it The ergodicity of a transformation $T$ with respect to a measure $\mu$ is equivalent to the following properties.

\bf 1. \it If $f$ is $\mu$-measurable and $Tf=f$, then $f$ is a constant almost everywhere.

\bf 2. \it If the measure $\nu \ll \mu$ and $\nu$ is invariant under $T$, then $\nu=c \,\mu$, $c\geq 0$.

\bf 3. \bf Average mixing: \it
for any measurable $A,B$
$$\frac 1 N \sum_1^N \mu(T^i A\cap B)\to \mu( A)\mu(B),$$
which is equivalent to the weak convergence $$\frac 1 N \sum_1^N T^i \to_w \Theta.$$
\rm

\bf 4. \bf Von Neumann's theorem \rm on strong convergence of ergodic averages: \rm $$\frac 1 N \sum_1^N T^i \to_s \Theta.$$

\bf 5. \bf Birkhoff's ergodic theorem \rm (convergence almost everywhere).

\it Let $f\in L_1$ and $T$ be an ergodic transformation.
Then for almost all $x$ we have the convergence
$$\frac 1 N \sum_1^N f(T^ix) \to \int f d\mu.$$
(\it Professionals say so: for an ergodic dynamical system the temporal averages converge to the spatial average.)
\rm

\vspace {4mm}
(1) The first assertion follows from the fact that $\mu(\{x: f(x)>c\})$ is 0 or 1 for all $c$,
hence $f=const= \inf$ over all $a$ for which $\mu(\{x: f(x)>a\})=0$.

(2) The second assertion. The Radon-Nikodym derivative (in other words, the density) of the measure $\nu$ with respect to $\mu$
is invariant with respect to $T$, therefore, it is constant (here we applied the first assertion).
We obtain $\nu=c \mu$, where $c$ is a non-negative number.

(3) Let us prove the average  mixing. Let us fix a measurable set $B$. From any sequence, one can select a subsequence $N_m$ such that for any $A$ the convergence
$$s_m(A)=\frac 1 {N_m} \sum_{i=1}^{N_m} \mu(T^i A\cap B)\to \nu(A),$$
where $\nu$ is some measure (by the way, why is that?).

By averaging $\nu(A)=\nu(TA)$, indeed,
$$ |s_m(TA)-s_m(A)|\leq \frac 2 {N_m},$$
therefore $s_m(TA)\to \nu(A)$ and $s_m(TA)\to \nu(TA)$. Thus, the measure $\nu$ is invariant with respect to $T$.

The measure $\nu$ is absolutely continuous with respect to $\mu$. This is evident from the equality
$$\mu=\nu+\nu',$$ where $\nu'$ is defined in the same way as $\nu$, but instead of $B$ we consider $B'=X\setminus B$.

Thus, by virtue of assertion 2, we obtain $\nu= const\, \mu = \mu(B) \mu$. The second equality is obtained if we substitute $X$ for
$A$. We have
$$\frac 1 {N_m} \sum_{i=1}^{N_m} \mu(T^i A\cap B)\to \mu(A)\mu(B),$$
and since the sequence $N_m$ was chosen as a subsequence of an arbitrary sequence,
we have
$$\frac 1 {N} \sum_{i=1}^{N} \mu(T^i A\cap B)\to \mu(A)\mu(B), \ \ N\to\infty.$$

We present an operator proof of a similar fact: If $T$ is an ergodic transformation, then
$$\frac 1 N \sum_1^N T^i \to_w \Theta.$$
We use the fact that any norm-bounded family of operators in $L_2$ has a limit point
in the weak operator topology.
Let $\frac 1 {N_m} \sum_1^{N_m} T^i \to_w P$. For any $f\in L_2$, in view of averaging, we have $TPf=Pf$,
therefore, by virtue of
the first assertion, $Pf=const=\int f d\mu$ (prove the second equality). This means that $P=\Theta$.

(4) We prove the strong convergence of $$\frac 1 N \sum_1^N T^i \to_s 0,$$
restricting ourselves to the space of functions with zero mean.

The strong convergence of $P_j\to_s 0$ means that $P_j^\ast P_j\to_w 0$. Indeed,
$$\|P_jf\|^2=(P_jf,P_jf)=(P_j^\ast P_jf,f)\to 0, \ \ \|P_jf\|\to 0.$$

For $P_N=\frac 1 N \sum_1^N T^i$ we have

$$ \|TP_N^\ast P_N - P_N^\ast P_N\|\to 0,$$
therefore the limit operator $Q$ for the sequence $P_{N_i}$
satisfies the condition
$$ TQ=Q.$$
Due to the ergodicity of $T$ we obtain $Q=0$, for the function $f$ with zero mean we have $$Qf=TQf=const=0.$$
We obtain $$P_N^\ast P_N \to 0, \ \ \|P_Nf\|\to 0,$$ as required.

\section{Birkhoff's Theorem}
\bf Theorem 3.1.
\it Let $f\in L_1$ and $T$ be an ergodic automorphism of a probability space.
Then for almost all $x$ we have the convergence
$$\frac 1 N \sum_1^N f(T^ix) \to \int f d\mu.$$ \rm

\vspace{3mm}
It suffices to consider the case
$\int f d\mu=0$, since the general case is trivially reducible to it. Let
$$X_a=\left\{x: \limsup_N \frac 1 N \sum_1^N f(T^ix) > a\right\}.$$
Obviously, $TX_a=X_a$, therefore, due to the ergodicity of $T$, the measure of $X_a$ is equal to 0 or 1.
Let us prove that the first case is impossible for $a>0$. Let it not be so. Then for (almost) every $x$
there is a minimal natural number $N(x)$ such that
$$ \ \frac 1 {N(x)} \sum_0^{N(x)-1} f(T^ix) > a. \eqno (\ast)$$
The function $N(x)$ is measurable, therefore, recalling real analysis, for any ${\varepsilon}>0$ we find $N$ such that
$$\mu(Y_N)> 1-{\varepsilon}, \ \ Y_N=\{x: N(x)<N\}.$$

Since we know the Rokhlin-Halmos lemma, it will be convenient for us to use it. We recall it. 

\vspace{2mm}
\it For an ergodic transformation $T$, ${\varepsilon}>0$ and a natural $H$, there exists a measurable set $B$ such that
the sets $T^iB$ are disjoint for $0\leq i\leq H-1$ and
$$\mu\left({\bigsqcup}_{i=0}^{H-1} T^iB\right)>1-{\varepsilon}.$$\rm


\vspace{2mm}
Take $H\gg  N \gg  1$ and the set $B$ from the lemma ($T^iB$ are disjoint for $0\leq i\leq H$).

Represent the phase space $X$ as the union of the following sets (pieces of orbits),
$$ x, Tx, T^2x, \dots, T^{h(x)}x,$$ where $x$ runs through all $B$, and $h(x)>0$ is chosen so that
$T^{h(x)+1}x\in B$, but the points $Tx, T^2x, \dots, T^{h(x)}x$ do not belong to $B$. Note that such a representation
is guaranteed by the Lemma, and $h(x)\geq H$.

Now the most important thing begins. We walk along the set $ x, Tx, \dots, T^{h(x)}x,$ performing the following
procedure. If $x\in Y_N$, then $N(x)$ first points, starting with $x$, put into a bag. Look at the next point, if it does not belong to $Y_N$, skip it and look at the next point, etc. Having encountered a point $x'$ from $Y_N$,
we again put into the bag $ x', Tx', T^2x', \dots, T^{N(x')-1}x'$ and move on, guided by the described rules.
Recall that here $N(x')\leq N\ll H \leq h(x)$.
This activity stops when $x'=T^hx$ and $H-h <N$. Having performed the procedure for each point $x\in B$,
denoting the bag by $X'$, we note that the measure of $X'$ is greater than $ 1 - {\varepsilon} - \frac N H .$ That is, this measure
can be arbitrarily close to 1.

On all pieces $ x', Tx', T^2x', \dots, T^{N(x')-1}x'$, the average value
of the function $f$ is greater than $a>0$, if you forgot about this, see $(\ast)$. Therefore, the inequality
$$\int_{X'} f d\mu > a\mu(X')$$ holds. Recalling the real analysis, or more precisely the absolute continuity of the Lebesgue integral, we get a contradiction:
$$ 0=\int_{X} f d\mu > a/2>0. $$ $a>0$.

Now let us define for $a>0$
$$X_{-a}=\left\{x: \liminf_N \frac 1 N \sum_1^N f(T^ix) < -a\right\}.$$

Similarly (or considering $-f$ instead of $f$) we get $\mu(X_{-a})=0$.

We only have to agree with Birkhoff:

$$\mu\left(x: \lim_N \frac 1 N \sum_1^N f(T^ix) =0\right)=1.$$

\vspace{2mm}  Birkhoff's theorem significantly used in the proof of the following interesting statement.

\bf Krygin-Atkinson theorem. \it  Let  $f:X\to \Z$, \ $\int_X f d\mu=0$,
and $T$ be ergodic. Then for any    $N$ and almost all $x\in X$
there is $N(x)>N$ such that
$$\sum_{i=0}^{N(x)-1} f(T^ix) =0.$$ \rm

\newpage
\section{Properties equivalent to weak mixing}
Let us consider the following properties of the transformation $T$.

\vspace{2mm}
\bf 1. \it The transformation $T$ has weak mixing:
$T^{n_i}\to_w \Theta.$ 

\vspace{2mm}
\bf 2. \it The spectrum of $T$ in the space orthogonal to the constants is continuous (there are no eigenfunctions).

\vspace{2mm}
\bf 3. \it The product $T\times T$ is ergodic with respect to $\mu\times\mu$.

\vspace{2mm}
\bf 4. \it The transformation $T$ has almost mixing:
$$\forall A,B\in{\mathcal B} \ \ \frac 1 {N} \sum_{i=1}^{N} |\mu(T^i A\cap B)- \mu(A)\mu(B)|\to 0, \ \ N\to\infty.$$
We can say this: most powers of $T$ are close to $\Theta$.

\vspace{2mm}
\bf 5. \it The transformation $T$ has the property: for any ergodic $S$, the product $S\times T$ is ergodic.

\vspace{2mm} 
\bf Theorem 4.1. \it Properties 1-5 are equivalent. \rm

\vspace{2mm}
Proof. The implication 1$\to$ 2 is obvious.

Show 2$\to$3 (not 3 implies not 2). Let $T$ be an ergodic transformation and the product
$T\times T$ is not ergodic.

Then there exists a set $D\subset X\times X$ invariant under $T\times T$, and
$0< \mu\times\mu (C) = c <1$. Let $K(x,y)=\chi_C(x,y)$. Consider the operator $P:L_2(\mu)\to L_2(\mu)$,
defined by the formula $$ Pf(y)=\int_X K(x,y)f(x) d\mu.$$

We have $$P Const=P^\ast Const= c\, Const.$$
Indeed,
$$ P\,Const(y)=\int_X K(x,y)\,Const d\mu=Const \int_X K(x,y) d\mu,$$
moreover, $$g(y)= \int_X K(x,y) d\mu(x)$$ is invariant with respect to $T$,
therefore, due to the ergodicity of the transformation $T$, the function $g$ is  constant, obviously equal to $c$.  For the operator $P^\ast$ we have the same.

Thus, the operators $P,P^\ast$ map the space $H$ of functions with zero mean (it is orthogonal to the constants) to itself. Since $K(x,y)$ is not a constant, we obtain $P^\ast PH\neq \{0\}$. (Prove this by showing that otherwise $P=с\,\Theta$.)

The operator $P^\ast P$ is a compact self-adjoint operator commuting with the operator $T$. (The latter follows from the invariance of $K(Tx,Ty)=K(x,y)$.) By the Hilbert-Schmidt theorem, there exists an eigenfunction $v$:
$ P^\ast Pv= a v$, where the eigenvalue $a$ is nonzero. Since $T$ commutes with $P^\ast P$, all vectors
$T^i v$ are eigenfunctions of the operator $P^\ast P$ with eigenvalue $a$. Since
the operator $P^\ast P$ is compact, we obtain that the space $L$ generated by all vectors $T^i v$ is finite-dimensional.
We have $TL=L$, and from the course of linear algebra we know that the linear operator $T$ on the finite-dimensional complex
space $L$ has an eigenvector. Thus, we have established 2$\to$3.

We prove 3$\to$4. Let us denote $$ c_i=\mu(T^i A\cap B), \ \ c= \mu(A)\mu(B).$$
From the ergodicity of $T$ as $N\to\infty$ we have
$$\frac 1 {N} \sum_{i=1}^{N} c_i\to c, $$
From the ergodicity of $T\times T$ --
$$\frac 1 {N} \sum_{i=1}^{N} c_i^2\to c^2. $$
(why?). But from what has been said it follows that
$$\frac 1 {N} \sum_{i=1}^{N} (c_i - c)^2\to 0, \ N\to \infty,$$
whence, due to the boundedness of $c_i$, we have
$$\frac 1 {N} \sum_{i=1}^{N} |c_i - c|\to 0, \ N\to \infty,$$
which is what was required.

The implication 4$\to$1 seems informally obvious. Exercise.

The implication 5$\to$3 is trivial.

4$\to$5. The idea of the proof is as follows: most of $T^i$ are close to $\Theta$, the operators
$\frac 1 {N} \sum_{i=1}^{N} S^i$ tend to $\Theta$ (ergodicity of $S$), therefore,
$$\frac 1 {N} \sum_{i=1}^{N} S^i\otimes T^i \to_w \Theta\otimes\Theta,$$
and this (mixing on average on cylinders)
is equivalent to the ergodicity of $S\times T$ with respect to $\mu\times\mu$.

\section{ On typical properties of transformations}
In functional analysis, the Baire theorem on categories plays a significant role (in the proof of the Banach-Steinhaus theorem and the Banach theorem on the inverse operator).
Baire categories have found applications in ergodic theory.
All transformations form a group $ Aut$, on which the complete Halmos metric is defined:
$$\rho (S,T)=\sum_{i=1}^{\infty}  2^{-1}\left[\mu(SB_i\Delta TB_i) +\mu(S^{-1}B_i\Delta T^{-1}B_i)\right],$$
where some fixed family of sets $\{B_i\}$ is dense in ${\mathcal B}$.
A family of transformations is typical (massive, generic) if its complement is a set of the first category, i.e. a countable union of nowhere dense sets. A property of a transformation is called typical if the set of all transformations with this property is typical. We will now list the properties that are important to us. The first two properties are typical.
Fix a standard probability space $(X,{\mathcal B},\mu)$ and consider
its automorphism group $ Aut $, equipped with the full Halmos metric $\rho $:
$$ \rho(S,T)=\sum_i 2^{-i}\left(\mu(SA_i\Delta TA_i)+\mu(S^{-1}A_i\Delta T^{-1}A_i)\right),$$
where the family of sets $\{A_i\}$ is dense in the algebra ${\mathcal B}$.
A property is said to be generic if the set of automorphisms (hereinafter referred to as transformations) with this property
contains some $G_\delta$-set dense in $Aut$.

Without proof, we use the fact that the metric is complete and the metric space $Aut$ is separable (see \cite{H}).
If $X$ is a segment $[0,1]$ with Lebesgue measure, then $Aut$ contains a dense family of all simplest rearrangements of segments.
They are obtained as follows: we split $X$ into segments of equal length and consider  all corresponding segments exchange transformations. We obtain a countable group of periodic transformations dense in $Aut$.

\vspace{3mm}
\bf 2.1. Typicalness and non-typicalness of  weak limits for powers.\rm

As Rokhlin and Halmos showed, the absence of mixing and the presence of weak mixing are typical properties. Let us see how to prove   more general facts.

A function of the operator $Q(T)$ is called admissible if it has the form 
$$Q(T)=a\Theta+\sum_i a_iT^i, \ \ a, a_i\geq 0, \ \sum_i a_i=1-a,$$
where $\Theta$ is the operator of orthoprojection onto the space of constants in $ L_2(X,{\mathcal B}, \mu)$.
Let us recall that a transformation $T$ satisfying the condition
$T^m\to\Theta$ as $m\to\infty$ is called mixing. We denote the transformation and the
transformation-induced operator in $ L_2(X,{\mathcal B}, \mu)$ in the same way.

\vspace{3mm}
\bf Theorem 5.1. \it For an infinite set $M\subset {\mathbf N}$ and admissible functions $Q$,$R$

(i) the set of transformations $T$ such that for some infinite subset $M(T)\subset M$
$T^m\to_w R(T)$ holds for $m\in M(T)$ , $m\to\infty$, is typical;

(ii) the set of transformations $T$ such that $ T^m\to_w Q(T)$ for $m\in M, m\to\infty$,
is a set of the first category.

\vspace{3mm}
\rm

For $R=I$, $Q=\Theta$, we obtain the above-mentioned results of Halmos and Rokhlin.

\rm

\vspace{3mm}
Proof (i). Fix a dense set of transformations $\{J_q\}$ in $Aut$, $q\in {\mathbf N}$. There exists an
infinite subset $M'\subset M$ and a weakly mixing transformation $S$ such that $S^m\to_w R(S)$,
$m\in M'$.

The desired $S$ can be realized as a rank-1 construction (see \S 10) with 
  $S^{h_j}\to_w R(S)$, $h_j\in M$.

Let $w$ denote the metric defining the weak operator topology. For any $n$ and $q$ we find a number $m=m(n,q)\in M$ and a neighborhood
$U(n,q)$ of the transformation $J_q^{-1}SJ_q$ such that the inequality
$$w(T^m,R(T))< \frac 1 n$$
is satisfied for all $T\in U(n,q)$. We obtain a $G_\delta$-set
$$W=\bigcap_n\bigcup_q U(n,q).$$ It is dense in $Aut$, since the conjugacy class
of the ergodic transformation $S$ is dense in $Aut$ (a consequence of the Rokhlin-Halmos lemma).

If $T\in W$, then for any $n$ there exists $q(n)$ such that for $m(n)=m(n,q(n))\in M$ the inequality $w(T^{m(n)},R(T))< \frac 1 n$ holds. Such $m(n)$ form a set $M(T)\subset M$.

We have thus obtained that $W$ consists of transformations satisfying the condition of item (i).

Assertion (ii) logically follows from (i). Indeed, fix an infinite $M$,
for a typical transformation $S$ there exists an infinite subset of $M$ on whose elements
$\{S: S^m\to_w R(S)\neq Q(S)\}$.
Thus, the condition $\{T: T^m\to_w Q(T)\}$ for $m\in M, m\to\infty$, is satisfied only for atypical $T$.

\newpage
\bf Theorem 5.2. \it Let $M(m)$ be a fixed sequence such that $M(m)\gg m$.
The following property of  transformations $T$ is typical: there is a  sequence $m_i\to\infty$ such that  any sequence $n_i$, $m_i\leq n_i\leq M(m_i)$, 
 is mixing:   $T^{n_i}\to_w \Theta.$ \rm

\vspace{3mm}
If, for example,  $$\LARGE M(m)=m!^{m!^{m!}},$$ then we see that  
sometimes our typical $T$  very ... very  long time  (from $m_i$ to some  $M>M(m_i)$)  must mix very well. 
However the typical transformation after such a  mixing epoch 
for an extremely long time becomes again rigid.

\vspace{3mm}
\bf Theorem 5.3. \it Let a sequence $r_n\gg n$ be fixed,
the following property of  transformations $T$ is typical:  for any $\eps>0$ there is $n$ for which   $\rho(Id, T^{nr})< \eps$ for all $r$, $1\leq r \leq r_n$.\rm

\vspace{2mm}
In fact there is a more general assertion.

\vspace{2mm}
\bf Theorem 5.4. \it Let $Q$ be admissible and  a sequence $N(p)\gg p$ be fixed.  The following property of $T$ is typical:  there is an infinite set $M\subset {\mathbf N}$ such that $T^{mp_m}\to_w Q(T)$ as $m\in M$,  $m\to\infty$, for any sequence $p_m$ satisfying  $1\leq p_m \leq N(p_m)$. \rm

\vspace{3mm} 
\bf Exercises. \rm -- Prove theorems 5.2, 5.3.
 
\vspace{2mm}
 -- Show that the weak limits of powers of the transformation form a semigroup.

\vspace{2mm} 
-- If the semigroup $W(T)$ of weak limits of powers of $T$ contains all operators of the form
$ \frac 1 2 (T^k + I)$, then the semigroup $W(T)$ contains all possible polynomials of the form $\sum_{i=1} a_iT^i$,
$a_i\geq 0$, $\sum_{i=1} a_i= 1$. 

\vspace{2mm}
-- If $T^{j^{2024}}\to I$, is it true that $T$ is a periodic transformation?

\newpage
\section{Lego to construct  transformations}
How to construct an ergodic transformation $T$,
having a weak limit of powers of the form $Q(T)=\sum_i a_i T^i$? Below is a  description of  constructions allowing to realize as weak limits of transformation powers  all admissible functions $Q$.

\vspace{2mm}
\bf Rank one constructions. \rm Let  
$$ \bar s_j=(s_j(1), s_j(2),\dots, s_j(r_j-1),s_j(r_j)), \ \ r_j>1, \ s_j(i)\geq 0,$$
a sequence of integer vectors, and  $h_1=1$. We will inductively define a transformation $T$  and its phase space.
At each step, what was defined at the previous ones does not change in the future.

At stage $j$, the  transformation $T$ cyclically 
permutes  
$$E_j, TE_j T^2,E_j,\dots, T^{h_j-1}E_j.$$
Such collection is called a tower.
For now $T$ is not defined on $T^{h_j-1}E_j.$
We cut the interval $E_j$ into $r_j$ intervals $E_j^1, E_j^2, E_j^3,\dots, E_j^{r_j}$ of the same measure.
Consider the columns
$$E_j^i, TE_j^i ,T^2 E_j^i,\dots, T^{h_j-1}E_j^i, \ i=1,2,\dots, r_j.$$
Then we construct $s_j(i)$ new intervals above the column with number $i$ and obtain a set of intervals
$$E_j^i, TE_j^i, T^2 E_j^i,\dots, T^{h_j+s_j(i)-1}E_j^i$$
(the intervals do not intersect and have the same measure).
For all $i<r_j$, we set
$$T^{h_j+s_j(i)}E_j^i = E_j^{i+1}.$$
Therefore, we have stacked the columns into a tower of stage $j+1$
$$E_{j+1}, TE_{j+1} T^2 E_{j+1},\dots, T^{h_{j+1}}E_{j+1},$$
where $$E_{j+1}= E^1_j,\
h_{j+1}=h_jr_j +\sum_{i=1}^{r_j}s_j(i).$$
Continuing the construction, we obtain a measure-preserving transformation $T$ on the union $X$ of all intervals.
If the measure $X$ is finite, we normalize it.

Mixing rank one constructions have the trivial centralizer (D. Ornstein, 
D. Rudolph, J. King; and V.V.  Ryzhikov, J.-P. Thouvenot for infinite mixing transformations).  S. Kalikow proved that (in the case of probability measure) mixing rank one constructions possess  2-fold mixing: for any measurable $A,B,C$
$$\mu(A\cap T^mB\cap T^{m+n}C)\to \mu(A)\mu(B)\mu(C), \ \ m,n\to +\infty.$$
By the way,  it is still unknown whether this is true in general --- \it Rokhlin's problem on multiple mixing remains open 75 years. \rm 

\vspace{2mm}
\bf Exercises. \rm  -- Prove that  rank one transformations are ergodic.

--Let $r_j=j$, $ s_j(i)=i$. Prove that $T^{h_j}\to\Theta$. Terry Adams proved
that this construction is mixing.  

-- Are there constructions such that $T^{n_j}\to\Theta$, but $T^{2n_j}\to I$?

-- Let $r_j=2j$, $s_j(i)=0$ as $1\leq i\leq j$ and  $s_j(i)=1$ as $j< i\leq 2j$.  This consruction called Katok's transformation  $T$.  Prove that
it is weakly mixing, non-mixing, and  weak closure of its powers have all operators in the form $\frac I 2 +\frac {T^k} 2$ (in fact all operatos $P=\sum_{i=1} a_iT^i$,
$a_i\geq 0$, $\sum_{i=1} a_i= 1$). 

\section{ Typical entropy invariants}
We define the following little modification of the Kirillov-Kushnirenko entropy.
Let $P=\{P_j\}$ be a sequence of finite subsets in a countable infinite group $G$.
For a measure-preserving action $T=\{T_g\}$ of $G$, we define the quantities
$$h_j(T,\xi)=\frac 1 {|P_j|} H\left(\bigvee_{p\in P_j}T_p\xi\right),$$
$$h_{P}(T,\xi)={\limsup_j} \ h_j(T,\xi),$$
$$h_{P}(T)=\sup_\xi h_{P}(T,\xi),$$
where $\xi$ denotes a finite measurable partition of $X$, $H(\xi)$ is the entropy of the partition $\xi$:
$$ H(\{C_1,C_2,\dots, C_n\})=-\sum_{i=1}^n \mu( C_i)\ln \mu( C_i).$$

We will be interested only in the case $|P_j|\to\infty$
(although the case of bounded cardinality also makes sense, directly related to the properties of the type  of multiple mixing).

\vspace{3mm}
\bf Entropy. \rm If $P_j=\{0,1,2,\dots,j\}$, then $h_P(T)$ is the classical
entropy $h(T)$, see \cite{KSF}.

\vspace{3mm}
\bf  $P$-entropy for large progressions. \rm Let's consider the following special case $G=Z$, when $P_j$ are increasing in size progressions:
$P_j=\{j,2j,\dots, L(j)j\}$, for some sequence $L(j)\to\infty$.

\vspace{3mm}
\bf Theorem 7.1. \it The set $\{S: h_P(S)=\infty\}$ is typical.\rm

Proof. Let $\{J_q\}$, $q\in \N$, be dense in $Aut$, and let
$T$ be the Bernoulli transform, denoted by $T_q=J_q^{-1}TJ_q$. The set $\{T_q\}$ is dense in $Aut$.
Fix a dense family $\{\xi_i\}$ of finite measurable partitions.

For any $n,q$ there exists $j=j(n,q)$ such that for all $i\leq n$ we have
$$ h_{j}(T_q,\xi_i)=\frac 1 {L_j} H(\bigvee_{n=1}^{L(j)} T^{nj}_q\xi_i) >H(\xi_i)-\frac 1 n. \eqno (\ast n\ast)$$
Indeed, $T_q$ is Bernoulli, we find a partition $\xi$ close to a fixed partition $\xi_i$, then
for some number $m(i,q,n,)$ the partitions $T^{nj}_q\xi$ are independent for all $n$ and $j>m(i,n,q)$.
This implies $(\ast n\ast)$ for all sufficiently large $j$.

We choose a neighborhood $U(n,q)$ of the transformation $T_q$ such that for all $S\in U(n,q)$ the inequality
$$ h_{j}(S,\xi_i) > H(\xi_i)-\frac 1 n.$$
The set
$$W=\bigcap_n\bigcup_q U(n,q),$$
is dense $G_\delta$. If $S\in W$, then for any $n$ there is $q(n)$
such that the inequality $$ h_{j(n,q(n))}(S,\xi_i)> H(\xi_i)-\frac 1 n$$ is satisfied for $i\leq n$.
And this leads to $h_P(S)=\infty$. Thus the set $\{S: h_P(S)=\infty\}$ is typical.

\vspace{3mm}
{ \bf Compact families with zero $P$-entropy.}
Denote by $E_0$ the class of transformations with zero entropy,
$K^{Aut}=\{J^{-1}SJ: S\in K, J\in Aut\}$.

\vspace{3mm}
\bf Theorem 7.2. \it If $K\subset E_0$ is a compact set and $(Aut,\rho)$, then the set $K^{Aut}$
is not typical.\rm

\vspace{3mm}
Proof. Fix a dense family of finite partitions $\xi_i$.
If $h(S)=0$, for any $i$ we have
$$h(S^j,\xi)=\lim_{L\to\infty}\frac 1 {L} H\left(\bigvee_{p=1}^L S^{jp}\xi_i\right)=0.$$
For $S\in K$ and $j$ we find a sequence of progressions $P(S)=\{P_{j}(S)\}$
$$ P_j(S)=\{j,2j,\dots, L_S(j)j\}$$
such that
$$\frac 1 {|P_j(S)|} H\left(\bigvee_{p\in P_j(S)}S^p\xi_i\right)<\frac 1 j$$
is satisfied for $i<j$.

Given the structure of the sets $P_j(S)$ and the fact that $K$ is compact, we find a sequence $L(j)\to\infty$
such that for any $S\in K$ and all sufficiently large $j$ we have $L(j)>L_S(j)$.
Then for a sequence $P$ of
expanding progressions $P_j=\{j,2j,\dots, L(j)j\}$ we have $h_P(S)=0$ for all $S\in K$ and thus for all
$h_P(T)=0$ for all $T\in K^{Aut}$. It follows from Theorem 7.1 that $K^{Aut}$ is a set of the first category,  which completes the proof.

\newpage
\section{Poisson suspensions with completely positive $P$-entropy}  
Ergodic theory studies large classes of dynamical systems that have an external origin. These include, in particular, Gaussian and Poisson automorphisms. The former are associated with the action of the group of all orthogonal operators on a space with a Gaussian measure, while the latter are the result of an injective embedding of the group of transformations preserving the sigma-finite measure into the group of transformations preserving the Poisson probability measure. For Gaussian automorphisms, see \cite{KSF}.

\bf The Poisson measure. \rm Consider the configuration space $X_\circ$, which consists of all infinite countable sets $x_\circ$ such that each above interval from  the spaces $X$ contains only a finite number of elements of the set $x_\circ$. 

The space $X_\circ$ is equipped with the Poisson  measure. We call its definition.  
For a  subset $A\subset X$ of a finite $\mu$-measure, we define configuration subsets
  $C(A,k)$, $k=0,1,2,\dots$,
  to $X_\circ$ by the formula
$$C(A,k)=\{x_\circ\in X_\circ \ : \ |x_\circ\cap A|=k\}.$$

 All possible finite intersections of the form $\cap_{i=1}^N C(A_i,k_i)$
    form a semiring.   
A Poisson measure $\mu_\circ$ is given on this semiring.
  Provided that the measurable sets $A_1, A_2,\dots, A_N$ do not intersect
and have a finite measure, we set
$$\mu_\circ\left(\bigcap_{i=1}^N C(A_i,k_i)\right)=\prod_{i=1}^N \frac {\mu(A_i)^{k_i}}{k_i!} e^{-\mu(A_i)}.\eqno (\circ)$$

The meaning of this formula is as follows: if the sets $A$, $B$ do not intersect, then
   probability $\mu_\circ(C(A,k))\cap C(B,m))$ of $k$ points of configuration $x_\circ$ in $A$ simultaneously
    and $m$ points of the configuration $x_\circ$ in $B$ is equal to the product of the probabilities $\mu_\circ(C(A,k))$ and 
$\mu_\circ( C(B,m))$.
   In other words, the events $C(A,k))$ and $C(B,m))$ are independent. Since the sets $A_1, A_2,\dots,A_N$ do not intersect,
 so the product appears in the formula $(\circ)$. Any element of the semiring is a finite union
semiring elements for which the Poisson measure is defined by $(\circ)$. The measure extends from the semiring to the Poisson configuration space 
$(X_\circ,\mu_\circ)$,
  isomorphic to the standard Lebesgue probability space.

An automorphism $T$ of the space $(X,\mu)$  naturally induces an automorphism 
$T_\circ$ of the space $(X_\circ,\mu_\circ)$,
  this $T_\circ$   is called  Poisson suspension.

\vspace{3mm}
\bf  Examples of $T_\circ$  with completely positive $P$-entropy. \rm
Back to rank one transformations, let $s_j(i)>L(j)h_j$,  $L(i)\to\infty$. 
Then $\mu (X_j)\to\infty$  and moreover for the correponding rank one construction $T$ the sets 
$$X_j, \ T^{h_j}X_j,\ T^{2h_j}X_j,\ \dots,\  T^{L(j)h_j}X_j$$ do not intersect.  The same  is automatically true for all  $A\subset X_j$.

Let  $C=C(A,k)$, where $A\subset X_{j_0}$, for the Poisson suspension $T_\circ$
we see that the sets
$$C, \ T^{h_j}_\circ C,\ T^{2h_j}_\circ C\  \dots, \ T^{L(j)h_j}_\circ C$$
are independent with respect to the Poisson measure.

Standard reasoning shows that $T_\circ$  has the completely positive $P$-entropy,
where $$P=\{ P_j \},  \ \ P_j=\{h_j, 2h_j,\dots, L(j)h_j\}.$$ 
From the above the following result follows.

\vspace{3mm}
\bf Theorem 8.1 \rm (V.V. Ryzhikov, J.-P. Thouvenot). \it For any ergodic 
transformation $S$ of zero  entropy, i.e. $h(S)=0$,  there is a Poisson suspension $T$
of zero  entropy  and a sequence of progression $P$ such that 
$h_P(S)=0$ and $h_P(T)=\infty$. \rm

\vspace{3mm}
\bf Remark.  \rm In this theorem Poisson suspensions with the same success can be replaced by Gaussian automorphisms (see \cite{KSF} for definitions).

\newpage
\section{Spectral theorem for unitary operators} Ivertible measure-preserving  transformations induce unitary operators
on $L_2$. All properties of the latters are completely determined by the so-called spectral measures. Let us consider  two  examples.

Let $ U:l_2({\mathbf Z})\to\l_2({\mathbf Z})$ be a shift in the space $l_2$ of two-sided sequences:
$$ Ue_n=e_{n+1}$$
(we can assume that an orthonormal basis $\{e_n\}$ is given in a Hilbert space, $U$ is a shift on this system,
which naturally extends to a unitary operator on the entire space).

Consider the operator: $V:L_2({\mathbf T},\sigma)\to L_2({\mathbf T},\sigma)$,
$$ Vf(z)=zf(z), |z|=1, $$
where ${\mathbf T}$ is the unit circle in the complex plane, $\sigma$ is the normalized Lebesgue measure on ${\mathbf T}$.
Note that the mapping $\Phi$,
$$ \Phi e_n=z^n,$$
(here $z^n$ is a function) realizes an isomorphism of the spaces $\l_2({\mathbf Z})$ and $L_2({\mathbf T},\sigma)$
and an isomorphism of the operators $U$ and $V$. It turns out that $V$ is the spectral representation of the operator $U$,
and $\sigma$ is its spectral measure.

\vspace{3mm}
\bf Remark. \rm One of Banach's open problems is formulated as follows: \it is there a transformation with simple Lebesgue spectrum,  i.e  a transformation that as an operator is isomorphic to the above operator  $U$?\rm

\vspace{3mm}
Let us  give now a very simple example of a unitary operator in one-dimensional space:  $Ux=-x$
and the multiplication operator $V$ in $L_2({\mathbf T},\sigma)$, given by the formula $ Vf(z)=zf(z),$
where $\sigma$ is a measure on ${\mathbf T}$ concentrated at a single point $-1$. This measure is the spectral measure
of the operator $U$ (and the operator $V$).

We now turn to a  general case.

\vspace{3mm}
\bf Theorem 9.1. \it If a unitary operator $U:H\to H$ has a cyclic vector, then it is isomorphic to the operator
$V:L_2({\mathbf T},\sigma)\to L_2({\mathbf T},\sigma)$, $ Vf(t)=tf(t), |t|=1 $, for some Borel measure on ${\mathbf T}$. \rm

Proof. Let $h$ be a cyclic vector, i.e. the closure of the space containing all $U^nh,\ n\in {\mathbf Z}$, is $H$. Let us define the function $\rho$ on ${\mathbf T}$ $$ \rho_N(z)=\frac 1 N \left(\sum_{i=0}^{N-1} z^iU^{-i}h,\ \sum_{j=0}^{N-1} z^jU^{-j}h\right) = \frac 1 N \left(\sum_{i,j=0}^{N- 1} z^{i-j}U^{j-i}h, \ h\right)\geq 0.$$ Since $$\int_{{\mathbf T}} z^{i-j} dm=0, \ i\neq j,\ \int_{{\mathbf T}} z^0 dm=1,$$ we obtain $$\int_{{\mathbf T}}\rho\ dm= \frac N N =1.$$
Consider the sequence of measures $\sigma_N$:
$$d\sigma_N=\rho_N(z) dm.$$
This sequence, like any sequence of normalized measures on a compact set, has a weak limit point, which is a normalized Borel measure on the circle ${\mathbf T}$:
$$\forall f\in C({\mathbf T}) \ \ \int_{{\mathbf T}} f d\sigma_{N_k}\to \int_{{\mathbf T}}f d\sigma.$$
The resulting measure is the desired one.

Since
$$\int_{{\mathbf T}} z^i d\sigma_{N_k}= \frac{N_k-|i|}{N_k}(U^ih,h)\to (U^ih,h)$$
and
$$\int_{{\mathbf T}} z^i d\sigma_{N_k} \to \int_{{\mathbf T}} z^i d\sigma,$$
we get
$$(U^ih,h)=\int_{{\mathbf T}}z^i d\sigma.$$
Thus,
$$(U^ih,U^j h)=\int_{{\mathbf T}}z^i \bar z^j d\sigma,$$
from which we arrive at the fact that the comparison
$$\Phi U^ih=z^i$$
preserves the scalar product and extends to a linear isometry of the spaces $H$ and $L_2(\sigma)$.

The resulting measure $\sigma$ is called the spectral measure, its definition depended on the choice of a cyclic vector,
however, all spectral measures for the operator $U$ are equivalent to each other.

Now we formulate the spectral theorem for the general case of unitary operators.

\vspace{3mm}
\bf Theorem 9.2. \it A unitary operator $U:H\to H$ on a separable Hilbert space
is isomorphic to the operator
$V$ acting in $L_2({\mathbf T}\times{\mathbf N},\sigma)$ for some Borel measure $\sigma$ on ${\mathbf T}\times{\mathbf N}$ by the formula
$$ Vf(z,n)=zf(z,n), \ \  f\in L_2({\mathbf T}\times{\mathbf N},\sigma).$$
\rm

The proof is easy to obtain from the previous theorem using
the decomposition of $H$ into an orthogonal sum of cyclic
subspaces.

\vspace{3mm}
\bf Exercises. \rm -- A unitary operator is completely determined (up to isomorphism) by some measure $\sigma$ on the unit circle ${\mathbf T}$ in the complex plane and the multiplicity function $m(z)$.
Let  a  normal  operator be the  multiplication by a function
$\phi: X\to {\mathbf C}$  on some measure space  $(X,\mu)$, how to find the corresponding spectral measure $\sigma$ on ${\mathbf C}$ and the multiplicity function  $m(z)$?

Solution.  Consider the graph of the function $\phi$ in $X\times {\mathbf C}$. Lift the measure $\mu$ onto it and project $\gamma$
on ${\mathbf C}$. We obtain $\sigma$. The measure $\gamma$ corresponds to a system of conditional measures $\gamma_z$.
The function $m(z)$ is the number of points on which the discrete measure $\gamma_z$ is concentrated. If the measure
$\gamma_z\neq 0$ is of another kind, we set $m(z)=\infty$. 

\vspace{2mm}
--  Decompose the space on which the unitary operator acts into the orthogonal sum of its cyclic  subspaces.

\vspace{2mm}
-- Find the spectral measure of the operator induced by rotating the circle by an angle $a$.

\vspace{2mm}
-- Prove von Neumann's theorem by using spectral theorem.

\section{Compact factors, Kronecker algebra}
We begin with some remarks on the eigenfunctions of the ergodic transformation of $T$.

\vspace{2mm}
1. \it If $Tf= \lambda f$, then $|f|=const$. \rm
\\
Since  $T|f|=|Tf|=|\lambda f|=|f|,$
and  $T$ is ergodic, we get $|f|=const$.

\vspace{2mm}
2. \it If $Tf= \lambda f$ and $Tg= \lambda g\neq 0$, then $\frac f g = \frac {Tf} {Tg} =Const$.
\\
In other words, the multiplicity of an eigenvalue of an ergodic transformation is equal to 1.

\vspace{2mm}
3. \it The eigenvalues of a measure-preserving transformation form a group.\rm
\\ This follows from the multiplicativity property of operators induced by a change of variable, for them the equality
$$T(fg)=Tf\, Tg.$$ holds. Then, if $Tf= \lambda_1 f$ and $Tg= \lambda_2 g$,
then $T(fg)= \lambda_1\lambda_2 fg.$

\vspace{2mm}
4. \it If $Sf= \lambda f\neq Const$ and $Tg= \lambda g\neq Const $, then the product $S\times T$ is not
ergodic. \rm
\\
The function $h(x,y)=\frac {f(x)} {g(y)} = \frac {Sf(x)} {Tg(y)} \neq Const$, therefore $S\times T$ is not ergodic.

\vspace{4mm}
\bf Compact factor coincides with  Kronecker algebra. \rm The smallest 
sigma-algebra ${\mathcal K}\subset {\mathcal B}$ with respect to which all eigenfunctions of the transformation $T$ are measurable is called the Kronecker algebra.

The compact factor ${\mathcal K}'$ of the transformation $T$ is the algebra of sets generated by
compact functions, i.e. functions $f\in L_2$ such that the orbit $\{T^i f\}$ is precompact in $L_2$.

\vspace{2mm}
\bf Theorem  10.1. \it The compact factor coincides with the Kronecker algebra:\ \ \ ${\mathcal K}={\mathcal K}'$.\rm

\vspace{2mm}
Proof outline. From the spectral theorem we obtain that the space on which the unitary operator induced by the transformation acts
is the orthogonal sum of two spaces:

the first is generated by the eigenvectors (corresponds to the discrete part of the spectrum),

and the second is the orthogonal complement 
to the first. 
The compact vectors form exactly the first space.

The space $L_2({\mathcal K}')$ consists of compact vectors, so it coincides with $L_2({\mathcal K})$, whence ${\mathcal K}={\mathcal K}'$.

 \section{ Progression  recurrence for  weakly mixing transformations}
The famous theorem of Szemeredi (he studied at the Faculty of Mechanics and Mathematics of Moscow State University) on progressions is equivalent to the fact that for any set $A$ of positive measure and an invertible measure-preserving transformation $T$ of the probability space, there is $i>0$ for which the property of multiple recurrence holds:
$$ \mu( A\cap T^iA\cap T^{2i}A\dots \cap T^{ki}A)\ >0.$$
This statement is easy to deduce from the theorem on progressions, but Furstenberg proved it using the methods of ergodic theory
and thus gave a completely new proof of Szemeredi's theorem.
The case $k=1$ is the Poincare recurrence theorem.
For $k=2$ and especially for $k>2$ the proof is not so easy as for  $k=1$.

\newpage
\bf Theorem 11.1. \it
For a weakly mixing transformation $T$ of the probability
space $(X,{\mathcal B},\mu)$ for any $ A,A_1,\dots A_k\in{\mathcal B}$ as $N\to\infty$ we have
$$ \frac 1 {N} \sum_{i=1}^{N} \mu( A\cap T^iA_1\cap T^{2i}A_2\dots \cap T^{ki}A_k)\ \to
\ \mu(A)\mu(A_1)\dots\mu(A_k). \eqno {\bf (a.p. Mix(k))}$$\rm

\vspace{3mm}
Let $k=2$. From any infinite subset of the natural numbers
we can always choose a subsequence $N(k)\to\infty$ so that for any $ A,A_1,A_2\in{\mathcal B}$
$$ \frac 1 {N(k)} \sum_{i=1}^{N(k)} \mu( A\cap T^iA_1\cap T^{2i}A_2)\ \to \ \nu(A\times A_1\times A_2),
$$
where $\nu(A\times A_1\times A_2)$ is for now a notation for the limit of expressions on the left. It is easy to see that
$\nu$ as a function on the cylinders $A\times A_1\times A_2$ that form a semiring is a measure.
The invariance is obvious
$$\nu(A\times A_1\times A_2) =\nu(TA\times TA_1\times TA_2),$$
and due to averaging there is an additional invariance
$$\nu(A\times A_1\times A_2) =\nu(A\times TA_1\times T^2A_2).$$
The measure $\nu$ as a measure on a cube is projected onto factors in the measure $\mu$. Such $T\times T\times T$-invariant measures with good projections are called self-joinings or joinings.

The weak mixing property is equivalent to the ergodicity of the transformation $T\times T^2$, which in turn means that any sets $ A_1\times A_2$ and $ B_1\times B_2$ are equidecomposable provided that $ \mu(A_1)\mu(A_2)=\mu(B_1)\mu(B_2).$
The latter means that
$$ A_1\times A_2\ ={\bigsqcup}_{n=1}^\infty (C_n\times D_n), $$
$$ B_1\times B_2\ ={\bigsqcup}_{n=1}^\infty (T^{p(n)}C_n\times T^{2p(n)}D_n),$$
where the equalities are satisfied up to sets of zero $\mu\times\mu$-measure.
We get
$$\sum_n \nu(A\times C_n\times D_n)=\sum_n \nu(A\times T^{p(n)}C_n\times T^{p(n)}D_n),$$
$$\nu(A\times A_1\times A_2) =\nu(A\times B_1\times B_2).$$
From the above it follows that $\nu=\mu\times\mu\times\mu$.
And this means
$$ \frac 1 {N} \sum_{i=1}^{N} \mu( A\cap T^iA_1\cap T^{2i}A_2)\ \to \ \mu(A)\mu(A_1)\mu(A_2).$$

\bf Exercise. \rm   Establish ${\bf (a.p. Mix(k))}$ by induction for all $k>2$.

\section {Double recurrence  for ergodic  transformations}

\bf Theorem (Roth, Furstenberg). \it For any set $A$ of positive measure
and an invertible measure-preserving transformation $T$ of the probability space $(X,{\mathcal B},\mu)$
there exists $i>0$ such that
$$ \mu( A\cap T^iA\cap T^{2i}A)\ >0.$$\rm

Proof(by J.-P. Thouvenot and  V.V. Ryzhikov, last century). Due to Rokhlin's theorem on the decomposition of an invariant measure into ergodic components, the general case reduces to the case of an ergodic transformation $T$, which we will now consider.
The ergodicity of the transformation $T$ is equivalent to 
$$\frac 1{N}\sum_{i=1}^{N} T^{i}f\to_w const =\int_X f\ d\mu, $$
for all $f\in L_\infty(\mu)$.

We define the operator $J:L_2(\mu)\to L_2(\mu\times\mu)$ by the formula
$$(Jf,g\otimes h)=\lim_{N_k}\frac 1{N_k} \sum_{i=1}^{N_k}\int_X f\ T^ig\ T^{2i}h\ d\mu,$$
where $f,g,h\in L_\infty(\mu)$.
The operator $J$ is well defined, since due to the ergodicity of the transformation $T$ we obtain
$$\lim_{N_k}\frac 1{N_k} \sum_{i=1}^{N_k}\int_X T^ig\ T^{2i}h\ d\mu =
\lim_{N_k}\frac 1{N_k}\sum_{i=1}^{N_k}\int_X g\ T^{i}h\ d\mu=\int_X g\ d\mu \int_X h\ d\mu.$$
The averaging that appears in the definition of the operator $J$ ensures the equality
$$(T\otimes T^2)J=J.$$

\bf Remark. \rm The equality that has appeared is directly related to the additional invariance discussed
$(Id\times T\times T^2)\nu=\nu$ for the self-joining $\nu$ corresponding to the operator $J$:
$$\nu(A\times B\times C) =(J\chi_A, \chi_B\otimes \chi_C).$$ Self-joinings and intertwining operators
in some cases are two sides of the same coin, the choice of sides is dictated by considerations of convenience.
The operator $J$ under consideration is the operator of conditional mathematical expectation:
$$(Jf)(x_1,x_2)=\int_X f(x)d\nu_{(x_1,x_2)}, $$
where $\nu_{(x_1,x_2)}$ is the family of conditional measures corresponding to the measure $\nu$.

Let us return to the proof. Since $(T\otimes T^2)Jf=Jf$, the image of the operator $J$ consists of functions that are fixed with respect to the operator $T\otimes T^2$. But the space of such fixed functions is a subspace of $V\otimes V$, where $V$ is generated by all eigenvectors of the operator $T$. This follows from the standard facts of the spectral theory of unitary operators, prove it using the spectral theorem.

Denote by $\pi$ the orthoprojection of $L_2(\mu)$ onto $V$.

\vspace{2mm}
\bf Exercise. \rm  The operator $\pi$ preserves the non-negativity of the function.

\vspace{2mm}
We have $$(T\otimes T^2)Jf=Jf, \ \ (\pi\otimes\pi)J=J,$$ $$ (Jf,g\otimes h)= (Jf,(\pi\otimes\pi)g\otimes h)=$$ $$\lim_{N_k}\frac 1{N_k} \sum_{i=1}^{N_k}\int_X f\ T^i\pi g\ T ^{2i}\pi h\ d\mu =\lim_{N_k}\frac 1{N_k} \sum_{i=1}^{N_k}\int_X \pi f\ R^i\pi g\ R^{2i}\pi h\ d\mu,$$ where $R$ is the constraint of $T$ on $V$.

\vspace{2mm}
\bf Exercise. \rm  We replaced $f$ with $\pi f$, why?

\vspace{2mm}
But the powers of the operator $R$ regularly turn out to be close to the identity
operator, in other words, (again the problem) the set of those $i$, when
$$\|R^i\pi f -\pi f\|<{\varepsilon}, \ \|R^{2i}\pi f -\pi f\|<{\varepsilon}$$
has positive density.

Setting $f=g=h=\chi_A$ subject to $\mu(A)>0$, we obtain
$$(J\chi_A,\chi_A\otimes \chi_A)=\lim_{N_k}\frac 1{N_k} \sum_{i=1}^{N_k}
\int_X \pi \chi_A\ R^i\pi \chi_A\ R^{2i}\pi\chi_A\ d\mu \ >\ 0.$$
We have thus shown that the inequality $ \mu( A\cap T^iA\cap T^{2i}A)\ >0$ holds
for an infinite set of values of $i$, which completes the proof.



\begin{thebibliography}{99}
\bibitem{H}  Halmos  P.R., Lectures on Ergodic Theory. Publ. Math.Soc. Japan, No. 3, Tokyo, 1956


\bibitem{KSF} Kornfeld I.P., Sinai Ya.G., Fomin S.V., Ergodic Theory, Moscow, 1980

\end{thebibliography}
\end{document}